\documentstyle[12pt]{article}
\title{\bf Asymptotic behavior of the smallest eigenvalue of matrices
associated with completely even functions (mod $r$)
\thanks{The first author was supported partially by National Science
Foundation of China Grant \#10971145 and by the Ph.D. Programs Foundation
of Ministry of Education of China Grant \#20100181110073 and by Program for New
Century Excellent Talents in University Grant \# NCET-06-0785 and by the
Lady Davis Fellowship at the Technion and the second author was
supported partially by the M. and M. Bank Mathematics Research Fund.}}
\author{Shaofang Hong \\
{\it Mathematical College, Sichuan University, Chengdu 610064,
P.R. China}\\
{\it E-mail address:} sfhong@scu.edu.cn; s-f.hong@tom.com; hongsf02@yahoo.com\\
Raphael Loewy\\
{\it Department of Mathematics, Technion-I.I.T., Haifa 32000, Israel}\\
{\it E-mail address:} loewy@techunix.technion.ac.il}
\date{}
\def\re{\par\hang\textindent}
\begin{document}
\maketitle
\newcommand{\SSS}{\stackrel}
\newcommand{\p}[2]{{\Phi_{#1}}(#2) }
\newcommand{\J}[2]{\left( \frac{#1}{#2}\right)}
\renewcommand{\thesection}{\arabic{section}}

\noindent{\bf Abstract}
In this paper we present systematically
analysis on the smallest eigenvalue of matrices associated with
completely even functions (mod $r$). We obtain several theorems on
the asymptotic behavior of the smallest eigenvalue of matrices
associated with completely even functions (mod $r$). In particular,
we get information on the asymptotic behavior of the smallest
eigenvalue of the famous Smith matrices. Finally some examples are
given to demonstrate the main results.

\noindent{\bf 1. Introduction and statements of results}\\

For any given arithmetical function $f$, we denote by $f(m, r)$ the
function $f$ evaluated at the greatest common divisor $(m, r)$ of
positive integers $m$ and $r$. Cohen [Co2] called the function $f(m,
r)$ a {\it completely even function} (mod $r$). Let $1\le x_1<
...<x_n< ...$ be a given arbitrary strictly increasing infinite
sequence of positive integers. For any integer $n\ge 1$, let
$S_n=\{x_1,...,x_n\}$. Let $I$ be the function defined for any
positive integer $m$ by $I(m):=m$. In 1876, Smith [S] published his
famous theorem showing that the determinant of the $n\times n$
matrix $[I(x_i, x_j)]$ on $S_n=\{1, ..., n\}$ is the product
$\prod^n_{k=1} \phi (k),$ where $\phi $ is Euler's totient function.
Smith also proved that if $S_n=\{1, ..., n\}$, then ${\rm
det}[f(i,j)]=\prod^n_{k=1} (f*\mu )(k)$, where $\mu $ is the
M\"obius function and $f*\mu $ is the Dirichlet convolution of $f$
and $\mu $. In 1972, Apostol [A2] extended Smith's result. In 1986,
McCarthy [Mc2] generalized Smith's and Apostol's results to the
class of even functions of $m$ (mod $r$), where $m$ and $r$ are
positive integers. A complex-valued function $\beta (m, r)$ is said
to be an {\it even function of} $m \ ({\rm mod} \ r)$ if $\beta (m,
r)=\beta ((m, r), r)$ for all values of $m$ [Co1, Co2]. Clearly a
completely even function (mod $r$) is an even function of $m \ ({\rm
mod} \ r)$, but the converse is not true. In 1993, Bourque and Ligh
[Bo-L1] extended the results of Smith, Apostol, and McCarthy. In
1999, Hong [Hon3] improved the lower bounds for the determinants of
matrices considered by Bourque and Ligh ([Bo-L1]). In 2002, Hong
[Hon4] generalized the results of Smith, Apostol, McCarthy and
Bourque and Ligh to certain classes of arithmetical functions.
Another kind of extension of Smith's determinant were obtained by
Codec\'a and Nair [CN] and Hilberdink [Hi].

Let $\varepsilon $ be a real number. Wintner [W] proved in 1944
that ${\rm lim\ sup}_{n\rightarrow
\infty}\Lambda_n(\varepsilon)<\infty$ if and only if
$\varepsilon>1$, where $\Lambda_n(\varepsilon)$ denotes the
largest eigenvalue of the matrix $N_n$ defined as follows:
$$
N_n:=\left({{(i, j)^{2\varepsilon}}\over {i^{\varepsilon}\cdot
j^{\varepsilon}}}\right)_{1\le i,j\le n}.
$$
Let $\lambda_n(\varepsilon)$ denote the smallest eigenvalue of the
matrix $N_n$. Lindqvist and Seip [LS] in 1998 use the work of
[He-L-S] about Riesz bases to investigate the asymptotic behavior of
$\lambda_n(\varepsilon)$ and $\Lambda_n(\varepsilon)$ as
$n\rightarrow\infty$. In particular, they got a sharp bound for
$\lambda_n(\varepsilon)$ and $\Lambda_n(\varepsilon)$. In 2004, Hong
and Loewy [Hon-Lo] made some progress in the study of asymptotic
behavior of the eigenvalues of the  $n\times n$ matrix $(\xi
_{\varepsilon }(x_i, x_j))$ on $S_n$, where $\xi _{\varepsilon }$ is
defined for any positive integer $m$ by $\xi _{\varepsilon
}(m):=m^{\varepsilon }$. It was proved in [Hon-Lo] that if
$0<\varepsilon\le 1$ and $q\ge 1$ is any fixed integer, then the
$q$-th smallest eigenvalue of the $n\times n$ matrix $(\xi
_{\varepsilon }(i, j))$ defined on the set $S_n=\{1, ..., n\}$
approaches zero as $n$ tends to infinity. Recently, Hong and Lee
[Hon-Le] studied the asymptotic behavior of the eigenvalues of the
reciprocal power LCM matrices and made some progress while Hong
[Hon11] got some results about asymptotic behavior of the largest
eigenvalue of matrices associated with completely even functions
(mod $r$). Notice also that Bhatia [Bh], Bhatia and Kosaki [Bh-K]
and Hong [Hon12] considered infinite divisibility of matrices
associated with multiplicative functions.

Given any set $S$ of positive integers, we define the class
$\tilde{\cal C}_S$ of arithmetical functions by
$$
\tilde {\cal C}_S:=\{f: (f*\mu )(d')>0 \ {\rm whenever} \ d'|x, \
{\rm for \ any } \ x\in S\}.
$$
For an arbitrary given strictly increasing infinite sequence
$\{x_i\}_{i=1}^{\infty}$ of positive integers, we define the class
$\tilde{\cal C}$ of arithmetical functions by
$$
\tilde{\cal C}:=\{f: (f*\mu )(d')>0 \ {\rm whenever} \ d'|x, \ {\rm
for \ any } \ x\in \{x_i\}_{i=1}^{\infty}\}.
$$
Let $S_n=\{x_1, ..., x_n \}$ for any integer $n\ge 1$. Then it is
clear that if $f\in \tilde{\cal C}$, then $f\in \tilde{\cal
C}_{S_n}$. In 1993, Bourque and Ligh ([Bo-L2]) showed that if $f\in
\tilde{\cal C}_{S_n}$, then the matrix $(f(x_i, x_j))$ (abbreviated
by $(f(S_n))$) is positive definite. Hong ([Hon1]) improved Bourque
and Ligh's bounds for det$(f(S_n))$ if $f\in \tilde{\cal C}_{S_n}$.
In [Hon7] and [Hon9], Hong obtained several results on the
nonsingularity of the matrix $(f(S_n))$. On the other hand, the
$n\times n$ matrix $(f[x_i, x_j])$ (abbreviated by $(f[S_n])$)
having $f$ evaluated at the least common multiple $[x_i,x_j]$ of
$x_i$ and $x_j$ as its $i,j$-entry on any set $S_n=\{x_1, ...,
x_n\}$ is not positive definite in general. It may even be singular.
In fact, Hong [Hon2] showed that for any integer $n\ge 8$, there
exists a gcd-closed set $S_n=\{x_1, ..., x_n\}$ (i.e. $(x_i, x_j)\in
S_n$ for all $1\le i, j\le n$) such that the $n\times n$ matrix
$(I[S_n])$ on $S_n$ is singular. It should be remarked that Cao
[Ca], Hong [Hon6, Hon8] and Hong, Shum and Sun [Hon-S-S] provided
several results on the nonsingularity of the $n\times n$ matrix
$(\xi _{\varepsilon }[S_n])$, where $\varepsilon $ is a positive
integer. We note also that Li [L] and Hong and Lee [Hon-Le] gave
partial answers to Hong's conjecture [Hon7] of real number power LCM
matrices. From Bourque and Ligh's result [Bo-L3] we can see that if
$S_n$ is a factor-closed set (i.e. it contains every divisor of $x$
for any $x\in S_n$) and $f$ is a multiplicative function such that
$(f*\mu )(d')$ is a nonzero integer whenever $d'|{\rm lcm}(S_n)$,
then the matrix $(f(x_i, x_j))$ divides the matrix $(f[x_i, x_j])$
in the ring $M_n({\bf Z})$ of $n\times n$ matrices over the
integers. Note also that Hong [Hon5] showed that for any {\it
multiple-closed set} $S_n$ (i.e. $y\in S_n$ whenever $x|y|{\rm
lcm}(S_n)$ for any $x\in S_n$, where ${\rm lcm} (S_n)$ means the
least common multiple of all elements in $S_n$) and for any {\it
divisor chain} $S_n$ (i.e. $x_1|...|x_n$), if $f$ is a completely
multiplicative function such that $(f*\mu )(d')\in {\bf
Z}\backslash\{0\}$ whenever $d'|{\rm lcm}(S_n)$, then the matrix
$(f(x_i, x_j))$ divides the matrix $(f[x_i, x_j])$ in the ring
$M_n({\bf Z})$. But such a factorization is no longer true if $f$ is
multiplicative. Some other factorization theorems about power GCD
matrices and power LCM matrices are obtained by Hong [Hon10], by
Haukkanen and Korkee [HK], by Hong, Zhao and Yin [HonZY], by Feng,
Hong and Zhao [FHZ], by Tan [T],  by Tan and Lin [TL], by Tan, Lin
and Liu [TLL] and by Xu and Li [XL].

For any given set $S$ of positive integers, it is natural to
consider the following class of arithmetical functions:
$$
{\cal C}_S:=\{f: (f*\mu )(d')\ge 0 \ {\rm whenever} \ d'|x, \ {\rm
for \ any } \ x\in S\}.
$$
In the meantime, associated with an arbitrary given strictly
increasing infinite sequence $\{x_i\}_{i=1}^{\infty}$ of positive
integers we define the following natural class of arithmetical
functions:
$$
{\cal C}:=\{f: (f*\mu )(d')\ge 0 \ {\rm whenever} \ d'|x, \ {\rm
for \ any } \ x\in \{x_i\}_{i=1}^{\infty}\}.
$$
Then it is easy to see that if $f\in {\cal C}$, then $f\in {\cal
C}_{S_n}$. Clearly $\tilde{\cal C}_S\subset {\cal C}_S$ for any
given set $S$ of positive integers, and $\tilde{\cal C}\subset
{\cal C}$ for any given strictly increasing infinite sequence
$\{x_i\}_{i=1}^{\infty}$ of positive integers. Obviously for any
given set $S$ of positive integers, $\tilde{\cal C}_S$ and ${\cal
C}_S$ are closed under addition and with respect to Dirichlet
convolution, and for any given strictly increasing infinite
sequence $\{x_i\}_{i=1}^{\infty}$ of positive integers,
$\tilde{\cal C}$ and ${\cal C}$ are closed under addition and with
respect to Dirichlet convolution. Note that $\mu \not\in {\cal
C}_S$ for any given set $S$ of positive integers containing at
least one prime, and $\mu \not\in {\cal C}$ for any given strictly
increasing infinite sequence $\{x_i\}_{i=1}^{\infty}$ of positive
integers containing at least one prime. However, we have the
following result (Theorem 1.1 below).

Let $c\ge 0$ be an integer. For any arithmetical function $f$,
define its $c$-th Dirichlet convolution, denoted by $f^{(c)}$,
inductively as follows:

$f^{(0)}:=\delta$ and $f^{(c)}:=f^{(c-1)}*f$ if $c\ge
1$,\\
where $\delta$ is the function defined for any positive integer
$m$ by
$$\delta(m):=\left\{\begin{array}{rl}
&1, \ \ {\rm if} \  m=1;\\
&0, \ \ {\rm otherwise}.
\end{array}\right.$$
Note that $f*\delta =f$ for any arithmetical function $f$ and
$$f^{(c)}:=\underbrace{f*...*f}_{c\ {\rm times}}.$$

For any integer $c\ge 1$, let
$$
{\bf Z}_{>0}^c:=\{(x_1, ..., x_c): 0<x_i\in {\bf Z}, \ {\rm for} \
i=1, ..., c\}.
$$

\noindent{\bf Theorem 1.1.} {\it Let $c\ge 1$ and $d\ge 0$ be
integers. If $f_1, ..., f_c$ are distinct arithmetical functions and
$(l_1,...,l_c)\in {\bf Z}_{>0}^c$, then each of the following is
true:}

(i). {\it Let $\{x_i\}_{i=1}^{\infty }$ be any given strictly
increasing infinite sequence of positive integers. If $f_1, ...,
f_c\in {\cal C}_{S_n}$ (resp. $f_1, ..., f_c\in {\cal C}$) and
$l_1+...+l_c>d$, then we have $f_1^{(l_1)}*...*f_c^{(l_c)}*\mu
^{(d)}\in {\cal C}_{S_n}$ (resp. $f_1^{(l_1)}*...*f_c^{(l_c)}*\mu
^{(d)}\in {\cal C}$);}

(ii). {\it For any prime $p$, we have
$$(f_1^{(l_1)}*...*f_c^{(l_c)}*\mu ^{(d)})(p)=\sum_{i=1}^cl_if_i(p)f_i(1)^{l_i-1}
\prod_{\SSS{j=1}{j\ne i}}^cf_j(1)^{l_j}
-d\prod_{i=1}^cf_i(1)^{l_i}.
$$
Furthermore, if $f_1, ..., f_c$ are multiplicative, then we have}
$$(f_1^{(l_1)}*...*f_c^{(l_c)}*\mu ^{(d)})(p)=\left\{\begin{array}
{cc}\sum_{i=1}^cl_if_i(p)-d
& \ {\it if} \ f_i(1)=1 \ {\it for \ all} \ 1\le i\le c,\\
0 & \ {\it if} \ f_i(1)=0 \ {\it for \ some} \ 1\le i\le c.
\end{array}\right.$$ \\

We remark that if the condition $l_1+...+l_c>d$ is suppressed,
then Theorem 1.1 (i) fails to be true. For example, let
$c=l_1=d=1$. Take $f_1=\phi$. Then $\phi \in {\cal C}_S$ for any
given set $S$ of positive integers and $\phi \in {\cal C}$ for any
given strictly increasing infinite sequence
$\{x_i\}_{i=1}^{\infty}$ of positive integers. But we have
$f_1*\mu=\phi *\mu \not\in {\cal C}_S$ for any given set $S$ of
positive integers containing at least one even number and $f_1*\mu
\not\in {\cal C}$ for any given strictly increasing infinite
sequence $\{x_i\}_{i=1}^{\infty}$ of positive integers containing
at least one even number because $(\phi *\mu ^{(2)})(2)=-1$.

Using Theorem 1.1 as well as [Theorem 1, Hon1] and by a
continuity argument, we can prove the following result.\\

\noindent{\bf Theorem 1.2.} {\it Let $c\ge 1$ and $d\ge 0$ be
integers and $S_n=\{x_1, ..., x_n\}$ be a set of $n$ distinct
positive integers. If $f_1, ..., f_c\in {\cal C}_{S_n}$ are distinct
and $(l_1,...,l_c)\in {\bf Z}_{>0}^c$ satisfies $l_1+...+l_c>d$,
then each of the following is true:}

(i). $\displaystyle \prod_{k=1}^n\sum_{\SSS{d'|x_k} {d'\not|x_t, \
x_t<x_k}}(f_1^{(l_1)}*...*f_c^{(l_c)}*\mu ^{(d+1)})(d')\le {\rm
det}((f_1^{(l_1)}*...*f_c^{(l_c)}*\mu ^{(d)})(x_i, x_j))\le
\prod_{k=1}^n(f_1^{(l_1)}*...*f_c^{(l_c)}*\mu ^{(d)})(x_k)$;

(ii). {\it The $n\times n$ matrix
$((f_1^{(l_1)}*...*f_c^{(l_c)}*\mu ^{(d)})(x_i, x_j))$
is positive semi-definite.}\\

Now let $\{x_i\}_{i=1}^{\infty}$ be an arbitrary given strictly
increasing infinite sequence of positive integers and let
$S_n=\{x_1, ..., x_n \}$ for any integer $n\ge 1$. Let $1\le q\le
n$ be a fixed integer and $c\ge 1$ and $d\ge 0$ be integers. Let
$(l_1,...,l_c)\in {\bf Z}_{>0}^c$ satisfy $l_1+...+l_c>d$ and
$f_1, ..., f_c\in {\cal C}$ be distinct. In the present paper, we
investigate the asymptotic behavior of the $q$-th smallest
eigenvalue of the matrix $((f_1^{(l_1)}*...*f_c^{(l_c)}*\mu
^{(d)})(S_n))$. Let $\lambda _n^{(1)}(l_1, ..., l_c, d)\le ...\le
\lambda _n^{(n)}(l_1, ..., l_c, d)$ be the eigenvalues of the
matrix $((f_1^{(l_1)}*...*f_c^{(l_c)}*\mu ^{(d)})(x_i, x_j))$
defined on the set $S_n=\{x_1, ..., x_n\}$. By Theorem 1.2 (ii) we
have
$$
\lambda _n^{(q)}(l_1, ..., l_c, d)\ge 0.
$$
But by Cauchy's interlacing inequalities (see [Hor-J1] and a
new proof of it, see [Hw]) we have
$$
\lambda _{n+1}^{(q)}(l_1, ..., l_c, d)\le \lambda _n^{(q)}(l_1,
..., l_c, d).
$$
Thus the sequence $\{\lambda _n^{(q)}(l_1, ...,
l_c, d)\}_{n=q}^{\infty }$ is a non-increasing infinite sequence
of nonnegative real numbers
and so it is convergent. Namely, we have \\

\noindent{\bf Theorem 1.3.} {\it Let $\{x_i\}_{i=1}^{\infty}$ be an
arbitrary given strictly increasing infinite sequence of positive
integers. Let $c\ge 1$ and $d\ge 0$ be integers and $q\ge 1$ be a
given arbitrary integer.  Let $f_1, ..., f_c\in {\cal C}$ be
distinct and $(l_1,...,l_c)\in {\bf Z}_{>0}^c$ satisfy
$l_1+...+l_c>d$. Let ${\lambda}_n^{(1)}(l_1, ..., l_c, d)\le ...\le
\lambda _n^{(n)}(l_1, ..., l_c, d)$ be the eigenvalues of the
$n\times n$ matrix $((f_1^{(l_1)}*...*f_c^{(l_c)}*\mu ^{(d)})(x_i,
x_j))$ defined on the set $S_n=\{x_1, ..., x_n\}$. Then the sequence
$\{\lambda _n^{(q)}(l_1, ..., l_c, d)\}_{n=q}^{\infty }$ converges
and}
$$
{\rm lim}_{n\rightarrow \infty }\lambda _n^{(q)}(l_1, ..., l_c,
d)\ge 0.
$$

Let $\{y_i\}_{i=1}^{\infty }$ be a strictly increasing infinite
sequence of positive integers. We say that $f$ is {\it increasing
on the sequence $\{y_i\}_{i=1}^{\infty }$} if $f(y_i)\le f(y_j)$
whenever $1\le i<j$. For an arbitrary strictly increasing infinite
sequence $\{x_i\}_{i=1}^{\infty }$ of positive integers satisfying
that $(x_i, x_j)=x$ for any $i\ne j$, where $x\ge 1$ is an
integer, we have the
following result.\\

\noindent{\bf Theorem 1.4.} {\it Let $x$ be a positive integer and
$\{x_i\}_{i=1}^{\infty}$ be a strictly increasing infinite sequence
of positive integers satisfying that for every $i\ne j, \ (x_i,
x_j)=x$. Assume that $f\in {\cal C}$ and is increasing on the
sequence $\{x_i\}_{i=1}^{\infty}$. Let ${\lambda}_n^{(1)}$ be the
smallest eigenvalue of the $n\times n$ matrix $(f(x_i, x_j))$
defined on the set $\{x_1, ..., x_n\}$. Then each of the following
holds:}

(i). {\it If $f(x)=0$, or $f(x)>0$ and $f(x_1)=f(x_2)$, then
${\lambda}_n^{(1)}=f(x_1)-f(x);$}

(ii). {\it If $f(x)>0$ and $f(x_1)<f(x_2)$, then}
$$f(x_1)-f(x)<{\lambda}_n^{(1)}<f(x_1)-f(x)+{{f(x)}\over
{1+\sum_{i=2}^{n}{{f(x)}\over {f(x_i)-f(x_1)}}}};$$

(iii). {\it If $f(x_1)=0$ then ${\lambda}_n^{(1)}=0$ for all $n\ge
1 $. If $f(x_1)>0$ and $\sum_{i=1}^{\infty }{1\over
{f(x_i)}}=\infty $, then we have ${\rm lim}_{n\rightarrow\infty
}{\lambda}_n^{(1)}=f(x_1)-f(x)$.}\\

From Theorem 1.4 we can deduce the following result.\\

\noindent{\bf Theorem 1.5.} {\it Let $x$ be a positive integer. Let
$\{x_i\}_{i=1}^{\infty}$ be a strictly increasing infinite sequence
of positive integers satisfying the following conditions:}

(i). {\it For every $i\ne j, \ (x_i, x_j)=x$;}

(ii). {\it $\sum_{i=1}^{\infty}{1\over {x_i}}={\infty}$.\\
Let ${\lambda}_n^{(1)}$ be the smallest eigenvalue of the $n\times
n$ matrix $(f(x_i, x_j))$ defined on the set $S_n=\{x_1, ...,
x_n\}$. If $f\in {\cal C}$ and is increasing on the sequence
$\{x_i\}_{i=1}^{\infty}$ and $f(x_i)\le Cx_i$ for all $i\ge 1$,
where $C>0$ is a constant, then we have ${\rm
lim}_{n\rightarrow\infty }{\lambda}_n^{(1)}=f(x_1)-f(x)$.}\\

Let $b\ge 1$ be an integer. By the well-known Dirichlet's theorem
(see, for example, [A1] or [I-R]) there are infinitely many primes
in the arithmetic progression $\{1+bi\}_{i=0}^{\infty }$. In the
following let
$$
p_1(b)<...<p_n(b)<... \ \ \ \ \eqno (1-1)
$$
denote the primes in this arithmetic progression. Consequently,
for the arithmetic progression case, we have the following result.\\

\noindent{\bf Theorem 1.6.} {\it Let $a, b, c, q\ge 1$ and $d, e\ge
0$ be any given integers. Let $x_i=a+b(e+i-1)$ for $i\ge 1$. Let
$(l_1,...,l_c)\in {\bf Z}_{>0}^c$ satisfy $l_1+...+l_c>d$. Let $f_1,
..., f_c\in {\cal C}$ be distinct, multiplicative and increasing on
the sequence $\{p_i(b)\}_{i=1}^{\infty}$, where $p_i(b) \ (i\ge 1)$
is defined by (1-1). Let ${\lambda}_n^{(1)}(l_1, ..., l_c, d)\le
...\le \lambda _n^{(n)}(l_1, ..., l_c, d)$ be the eigenvalues of the
$n\times n$ matrix $((f_1^{(l_1)}*...*f_c^{(l_c)}*\mu ^{(d)})(a+bi,
a+bj))$ defined on the set $\{a+be, a+b(e+1), ..., a+b(e+n-1)\}$.}

(i). {\it If $(f_1^{(l_1)}*...*f_c^{(l_c)}*\mu ^{(d)})(p_i(b))=0$
for some $i\ge 1$, then for any large enough $n$ we have
${\lambda}_n^{(q)}(l_1, ..., l_c, d)=0$;}

(ii). {\it If $(f_1^{(l_1)}*...*f_c^{(l_c)}*\mu ^{(d)})(p_i(b))\ne
0$ for all $i\ge 1$ and $\sum_{i=1}^{\infty }{1\over
{f_1(p_i(b))+...+f_c(p_i(b))}}=\infty $, then we have ${\rm
lim}_{n\rightarrow\infty }{\lambda}_n^{(q)}(l_1, ..., l_c, d)=0$;}

(iii). {\it In particular, if for each $1\le j\le c$, there is a
positive constant $C_j$ such that $f_j(p_i(b))\le C_jp_i(b)$ for
all $i\ge 1$, then we have ${\rm lim}_{n\rightarrow\infty
}{\lambda}_n^{(q)}
(l_1, ..., l_c, d)=0$.}\\

Furthermore, applying again Cauchy's interlacing inequalities, it
follows from Theorems 1.3 and 1.6 that the following result holds.\\

\noindent{\bf Theorem 1.7.} {\it Let $a, b, c, q\ge 1$ and $d, e\ge
0$ be any given integers. Let $\{x_i\}_{i=1}^{\infty }$ be any given
strictly increasing infinite sequence of positive integers which
contains the arithmetic progression $\{a+bi\}_{i=e}^{\infty }$ as
its subsequence. Let $(l_1,...,l_c)\in {\bf Z}_{>0}^c$ satisfy
$l_1+...+l_c>d$. Let $f_1, ..., f_c\in {\cal C}$ be distinct,
multiplicative and increasing on the sequence
$\{p_i(b)\}_{i=1}^{\infty}$, where $p_i(b) \ (i\ge 1)$ is defined by
(1-1). Let ${\lambda}_n^{(1)}(l_1, ..., l_c, d)\le ...\le \lambda
_n^{(n)}(l_1, ..., l_c, d)$ be the eigenvalues of the $n\times n$
matrix $((f_1^{(l_1)}*...*f_c^{(l_c)}*\mu ^{(d)})(x_i, x_j))$
defined on the set $S_n=\{x_1, ..., x_n\}$.}

(i). {\it If $(f_1^{(l_1)}*...*f_c^{(l_c)}*\mu ^{(d)})(p_i(b))=0$
for some $i\ge 1$, then for any large enough $n$ we have
${\lambda}_n^{(q)}(l_1, ..., l_c, d)=0$;}

(ii). {\it If $(f_1^{(l_1)}*...*f_c^{(l_c)}*\mu ^{(d)})(p_i(b))\ne
0$ for all $i\ge 1$ and $\sum_{i=1}^{\infty }{1\over
{f_1(p_i(b))+...+f_c(p_i(b))}}=\infty $, then we have ${\rm
lim}_{n\rightarrow\infty }{\lambda}_n^{(q)}(l_1, ..., l_c, d)=0$;}

(iii). {\it In particular, if for each $1\le j\le c$, there is a
positive constant $C_j$ such that $f_j(p_i(b))\le C_jp_i(b)$ for
all $i\ge 1$, then we have ${\rm lim}_{n\rightarrow\infty
}{\lambda}_n^{(q)}(l_1, ..., l_c, d)=0$.}\\

As a special case we have the following theorem.\\

\noindent{\bf Theorem 1.8.} {\it Let $a, b, c, q\ge 1$ and $d, e\ge
0$ be any given integers such that $c>d$. Let $\{x_i\}_{i=1}^{\infty
}$ be any given strictly increasing infinite sequence of positive
integers which contains the arithmetic progression
$\{a+bi\}_{i=e}^{\infty }$ as its subsequence. Let
${\lambda}_n^{(1)}(c,d)\le ...\le \lambda _n^{(n)}(c,d)$ be the
eigenvalues of the $n\times n$ matrix $((f^{(c)}*\mu ^{(d)})(x_i,
x_j))$ defined on the set $S_n=\{x_1, ..., x_n\}$. Let $f\in {\cal
C}$ be multiplicative and increasing on the sequence
$\{p_i(b)\}_{i=1}^{\infty}$, where $p_i(b) \ (i\ge 1)$ is defined by
(1-1).}

(i). {\it If $(f^{(c)}*\mu ^{(d)})(p_i(b))=0$ for some $i\ge 1$,
then for any large enough $n$ we have ${\lambda}_n^{(q)}(c,d)=0$;}

(ii). {\it If $(f^{(c)}*\mu ^{(d)})(p_i(b))\ne 0$ for all $i\ge 1$
and $\sum_{i=1}^{\infty }{1\over {f(p_i(b))}}=\infty $, then for
any given integer $q\ge 1$, we have ${\rm lim}_{n\rightarrow\infty
}{\lambda}_n^{(q)}(c,d)=0$;}

(iii). {\it In particular, if $f(p_i(b))\le Cp_i(b)$ for all $i\ge
1$, where $C>0$ is a constant, then for any given integer $q\ge
1$, we have ${\rm lim}_{n\rightarrow\infty }{\lambda}_n^{(q)}(c,d)=0$.}\\

\noindent{\bf Corollary 1.9.} {\it Let ${\lambda}_n^{(1)}\le ...\le
\lambda _n^{(n)}$ be the eigenvalues of the $n\times n$ matrix
$(f(i, j))$ defined on the set $S_n=\{1, ..., n\}$. If $f$ is an
increasing multiplicative function satisfying $(f*\mu )(y)\ge 0$ and
$f(y)\le Cy$ for all positive integers $y$, where $C>0$ is a
constant, then for any given integer $q\ge 1$, we have ${\rm
lim}_{n\rightarrow\infty }{\lambda}_n^{(q)}=0$.}\\

This paper is organized as follows. The details of the proofs of
Theorems 1.1-1.2 and 1.4-1.6 will be given in Section 2. In
Section 3 we give some examples to illustrate our results. The
final section is devoted to some open questions.

Throughout this paper, we let $E_n$ denote the $n\times n$ matrix
with all entries equal to 1. For the basic facts about arithmetical
functions, the readers are referred to [A1], [N] or [Mc1]. For a
comprehensive review of papers related to the matrices associated
with arithmetical functions not presented here, we refer to
[Hon-Le] and [Hon-Lo] as well as the papers listed there.\\

\noindent{\bf 2. The proofs of Theorems 1.1-1.2 and 1.4-1.6}\\

First we prove Theorem 1.1.\\

\noindent{\bf Proof of Theorem 1.1.} Clearly to prove Theorem 1.1 it
suffices to prove that for any prime $p$ and for any integer $l\ge
1$ and any (not necessarily distinct) arithmetical functions $g_1,
..., g_l$, we have
$$(g_1*...*g_l*\mu^{(d)})(p)=\sum_{i=1}^lg_1(1)...g_{i-1}(1)g_i(p)g_{i+1}(1)...g_l(1)
-dg_1(1)...g_l(1), \eqno (2-1)
$$
and if $g_1, ..., g_l\in {\cal C}_{S_n}$ (resp. $g_1, ..., g_l\in
{\cal C}$) and $l>d$, then we have
$$g_1*...*g_l*\mu ^{(d)}\in {\cal C}_{S_n} \ ({\rm resp.} \
g_1*...*g_l*\mu ^{(d)}\in {\cal C}). \eqno (2-2)$$
Furthermore, if $g_1, ..., g_l$ are multiplicative, then we have
$$(g_1*...*g_l*\mu^{(d)})(p)=\left\{\begin{array}{cc}\sum_{i=1}^lg_i(p)-d
& \ {\rm if} \ g_i(1)=1 \ {\rm for \ all} \ 1\le i\le l;\\
0 & \ {\rm if} \ g_i(1)=0 \ {\rm for \ some} \ 1\le i\le l.
\end{array}\right. \ \ \ \ (2-3)
$$

By the definition of Dirichlet convolution we have

$$\begin{array}{rl}
&(g_1*...*g_l*\mu^{(d)})(p)\\
=&\displaystyle \sum_{\SSS{r_1...r_l\bar r_1...\bar
r_d=p}{(r_1,...,r_l,\bar r_1,...,\bar r_d)\in {\bf Z}_{>0}^{l+d}}}
g_1(r_1)...g_l(r_l)\mu(\bar r_1)...\mu(\bar r_d)\\
=&\displaystyle\sum_{i=1}^lg_1(1)...g_{i-1}(1)g_i(p)g_{i+1}(1)...g_l(1)
\mu (1)^d+dg_1(1)...g_l(1)\mu (p)\mu (1)^{d-1}\\
=&\displaystyle\sum_{i=1}^lg_1(1)...g_{i-1}(1)g_i(p)g_{i+1}(1)...g_l(1)
-dg_1(1)...g_l(1).
\end{array}
$$
So (2-1) is proved. Further, if $f$ is multiplicative, then we
have $f(1)^2=f(1)$. So we have $f(1)=1$, or 0. Thus (2-3) follows
immediately.

Now consider (2-2). Since the proof for the case $g_1, ..., g_l\in
{\cal C}$ is completely similar to that of the case $g_1, ...,
g_l\in {\cal C}_{S_n}$, we only need to show (2-2) for the case
$g_1, ..., g_l\in {\cal C}_{S_n}$. In the following let $g_1, ...,
g_l\in {\cal C}_{S_n}$ and $l>d$. Now for any $x\in S_n$ and any
$r|x$, since $l\ge d+1$, we have
$$\begin{array}{rl}
&((g_1*...*g_l*\mu ^{(d)})*\mu )(r)\\
=&(g_1*...*g_l*\mu^{(d+1)})(r)\\
=&((g_1*\mu)*...(g_d*\mu )*(g_{d+1}*\mu )*g_{d+2}...*g_l)(r)\\
=&\displaystyle \sum_{\SSS{r_1...r_l=r}{(r_1,...,r_l)\in {\bf
Z}_{>0}^{l}}}(g_1*\mu)(r_1)...(g_{d+1}*\mu)(r_{d+1})g_{d+2}(r_{d+2})...g_l(r_l).
\ \ \ \ \ \ (2-4)
\end{array}
$$

For $1\le i\le d+1$, since $g_i\in {\cal C}_{S_n}$ and $r_i|x$, we
have
$$(g_i*\mu)(r_i)\ge 0. \ \ \ \ \eqno (2-5)$$
On the other hand, for $d+2\le j\le l$, $g_j\in {\cal C}_{S_n}$
together with $r_j|x$ implies that
$$
g_j(r_j)=\sum_{d'|r_j}(g_j*\mu )(d')\ge 0. \eqno (2-6)
$$
From (2-4)-(2-6) we then deduce that
$$
((g_1*...*g_l*\mu^{(d)})*\mu
)(r)\ge 0.
$$
Thus (2-2) holds. This completes the proof of Theorem 1.1. \ \ \ $\Box$\\

To prove Theorem 1.2 we need a result from [Hon1].\\

\noindent{\bf Lemma 2.1.} ([Theorem 1, Hon1]) {\it Let $S_n=\{x_1,
..., x_n\}$ be a set of $n$ distinct positive integers. If $g\in
\tilde{\cal C}_{S_n}$, then we have}

$${\rm det}(g(x_i, x_j))\ge \displaystyle
\prod_{k=1}^n\sum_{\SSS{d|x_k} {d\not|x_t, \ x_t<x_k}}(g*\mu
)(d).$$

We can now prove Theorem 1.2.\\

\noindent{\bf Proof of Theorem 1.2.} By Theorem 1.1 (i), to show
Theorem 1.2 we only need to show that if $f\in {\cal C}_{S_n}$, then
each of the following is true:

(i'). $\displaystyle \prod_{k=1}^n\sum_{\SSS{d'|x_k} {d'\not|x_t,
\ x_t<x_k}}(f*\mu)(d')\le {\rm det}(f(x_i, x_j))\le
\prod_{k=1}^nf(x_k)$;

(ii'). The $n\times n$ matrix $(f(x_i, x_j))$ is positive
semi-definite.

First we show the inequality on the left-hand side of (i'). Let
$f\in {\cal C}_{S_n}$. Pick $\epsilon
>0$ and $\bar f\in \tilde{\cal C}_{S_n}$. Then it is easy to see
that $f+\epsilon \bar f\in \tilde{\cal C}_{S_n}$. For an
arithmetical function $g$ and $1\le k\le n$, let

$$
\alpha _g(x_k):=\sum_{\SSS{d|x_k} {d\not|x_t, \ x_t<x_k}}(g*\mu
)(d).
$$

By Lemma 2.1 we have

$${\rm det}((f+\epsilon \bar f)(x_i, x_j))\ge \displaystyle
\prod_{k=1}^n\alpha _{f+\epsilon \bar f}(x_k). \ \ \ \eqno (2-7)$$

Note that both sides of (2-7) are polynomials in $\epsilon $.
Moreover, the constant coefficients of the left and right hand
sides are, respectively, ${\rm det}(f(x_i, x_j))$ and
$\prod_{k=1}^n\alpha _f(x_k)$. Since (2-7) holds for any
$\epsilon>0$, letting $\epsilon\rightarrow 0$ the left-hand
side of (i') is proved.

For any $1\le l\le n$, since $f\in {\cal C}_{S_n}$, then the
inequality on the left-hand side of (i') implies that the
determinant of any principal submatrix of $(f(x_i, x_j))$ is
nonnegative. This concludes part (ii'). From (ii') the inequality
on the right-hand side of (i') follows immediately. Hence
the proof of Theorem 1.2 is complete. \ \ \ $\Box$  \\

The following result is known.\\

\noindent{\bf Lemma 2.2.} {\it Let $n\ge 1$ be an integer and let
$a_1, ..., a_n\in R$, where $R$ is an arbitrary commutative ring.
Then we have}
$$
{\rm det}(E_n+{\rm diag}(a_1-1, ...,
a_n-1))=\prod_{i=1}^n(a_i-1)+\sum_{1\le i_1<...<i_{n-1}\le
n}\prod_{j=1}^{n-1}(a_{i_j}-1).
$$

In order to show Theorem 1.4 we need also the following lemma.\\

\noindent{\bf Lemma 2.3.} {\it Let $\{r_i\}_{i=1}^{\infty}$ be an
increasing infinite sequence of real numbers satisfying $r_1\ge 1$
and let ${\lambda}_n^{(1)}$ be the smallest eigenvalue of the
$n\times n$ matrix $E_n+{\rm diag}(r_1-1, ..., r_n-1)$. Then each of
the following holds:}

(i). {\it If $r_1=r_2$, then ${\lambda}_n^{(1)}=r_1-1.$}

(ii). {\it If $r_1<r_2$, then}
$$r_1-1<{\lambda}_n^{(1)}<r_1-1+{1\over
{1+\sum_{i=2}^{n}{1\over {r_i-r_1}}}}.$$

(iii). {\it If $\sum_{i=1}^{\infty}{1\over {r_i}}={\infty}$, then
${\rm lim}_{n\rightarrow\infty }{\lambda}_n^{(1)}=r_1-1$.}\\

\noindent{\bf Proof.} Clearly part (iii) follows immediately from
parts (i) and (ii). In what follows we show parts (i) and (ii).

Write
$$
F_n:=E_n+{\rm diag}(r_1-1, ..., r_n-1).
$$
Note that $F_n$ is positive semi-definite. Consider its
characteristic polynomial det$(\lambda I_n-F_n)$. By Lemma 2.2 we
have

$$\begin{array}{rl}
&(-1)^n{\rm det}(\lambda I_n-F_n)\\
=&{\rm det}(E_n+ {\rm
diag}(r_1-\lambda-1, ..., r_n-\lambda-1))\\
=&\displaystyle \prod_{i=1}^n(r_i-\lambda-1)+\sum_{1\le
i_1<...<i_{n-1}\le n}\prod_{j=1}^{n-1}(r_{i_j}-\lambda-1).
\end{array} \ \ \ \eqno (2-8)$$
We then deduce that if $\lambda <r_1-1$, then
$$(-1)^n{\rm
det}(\lambda I_n-F_n)>0
$$
and thus
$$
{\rm det}(\lambda I_n-F_n)\ne 0.
$$
So we have $\lambda _n^{(1)}\ge r_1-1$.

If $r_1=r_2$, then by (2-8) we have
$$
(\lambda-r_1+1)| {\rm det}(\lambda I_n-F_n). $$
It follows that
$\lambda _n^{(1)}=r_1-1$ and this concludes part (i).

Now let $r_2>r_1$. From (2-8) we deduce
$$
(-1)^n{\rm
det}((r_1-1)I_n-F_n)>0.
$$
This implies that $\lambda _n^{(1)}>r_1-1$. On the other hand, we
have
$$
F_n=(r_1-1)I_n+ E_n+{\rm diag}(0, r_2-r_1, ..., r_n-r_1).
$$

Let $\tilde\lambda _n^{(1)}$ be the smallest eigenvalue of the
$n\times n$ matrix $E_n+{\rm diag}(0, r_2-r_1, ..., r_n-r_1)$.
Then we have
$$
\lambda _n^{(1)}=r_1-1+\tilde\lambda _n^{(1)}. \eqno (2-9)
$$

Since $r_1<r_2$, the proofs of Lemma 2.2 and Corollary 2.3 of
[Hon-Lo] yield
$$\tilde\lambda _n^{(1)}<{1\over {1+\sum_{i=2}^n}{1\over {r_i-r_1}}}.\ \eqno (2-10)$$
So the right-hand side of the inequalities in part (ii) follows
immediately from (2-9) and (2-10). The proof of Lemma 2.3 is complete.\ \ \ \ $\Box$\\

We are now ready to prove Theorem 1.4.\\

\noindent{\bf Proof of Theorem 1.4.} By $f\in {\cal C}$ we have
$$f(x)=\sum_{d|x}(f*\mu
)(d)\ge 0.$$ If $f(x)=0$, then $(f(x_i, x_j))={\rm diag}(f(x_1),
..., f(x_n))$. Since $f$ is increasing on the sequence
$\{x_i\}_{i=1}^{\infty}$, we have $f(x_i)\ge f(x_1)$ for $1\le
i\le n$. Thus ${\lambda}_n^{(1)}=f(x_1)$. So Theorem 1.4 (i) is
true in this case. Now let $f(x)\ne 0$, so $f(x)>0$. Obviously we
have
$$
{1\over {f(x)}}(f(x_i, x_j))=E_n+{\rm diag}({{f(x_1)}\over
{f(x)}}-1, ..., {{f(x_n)}\over {f(x)}}-1). \ \ \ \  \ \eqno (2-11)
$$
For $1\le i\le n$, let $r_i={{f(x_i)}\over {f(x)}}$. Since $f$ is
increasing on the sequence $\{x_i\}_{i=1}^{\infty}$,
$\{r_i\}_{i=1}^{\infty}$ is an increasing infinite sequence of
real numbers. Since $x|x_1$ and $f\in {\cal C}$, we have
$$f(x_1)-f(x)=\sum_{d|x_1,\ d\not |x}(f*\mu )(d)\ge 0.$$
So $f(x_1)\ge f(x)$, namely, $r_1\ge 1$. Let $\bar
{\lambda}_n^{(1)}$ be the smallest eigenvalue of the $n\times n$
matrix ${1\over {f(x)}}(f(x_i, x_j))$ defined on the set
$S_n=\{x_1, ..., x_n\}$.

Suppose first $f(x_1)=f(x_2)$. Then $r_1=r_2$. Thus by (2-11) and
Lemma 2.3 (i) we have
$$\bar \lambda _n^{(1)}={{f(x_1)}\over {f(x)}}-1.$$
Theorem 1.4 (i) in this case then follows immediately from the
fact that
$$\lambda _n^{(1)} =f(x)\cdot \bar \lambda _n^{(1)}.  \ \ \eqno (2-12)$$
This completes the proof of Theorem 1.4 (i).

Let now $f(x_1)<f(x_2)$, i.e. $r_1<r_2$. By (2-11) and Lemma 2.3
(ii) we have
$${{f(x_1)}\over {f(x)}}-1<\bar \lambda _n^{(1)}<{{f(x_1)}\over {f(x)}}-1+{1\over
{1+\sum_{i=2}^{n}{{f(x)}\over {f(x_i)-f(x_1)}}}}.$$ So, by (2-12)
part (ii) of Theorem 1.4 follows.

Finally we show part (iii). If $f(x_1)=0$, then we have $f(x)=0$
because $f(x_1)\ge f(x)\ge 0$. Then by part (i) we have $\lambda
_n^{(1)}=0$ for $n\ge 1$. Thus Theorem 1.4 (iii) holds in this
case. If $f(x_1)>0$ and $\sum_{i=1}^{\infty }{1\over
{f(x_i)}}=\infty $, then part (iii) in this case follows
immediately from parts (i) and (ii). So part (iii) of Theorem 1.4
is proved. \ \ \ \ $\Box$\\

\noindent{\bf Proof of Theorem 1.5.} If $f(x_1)=0$, then by Theorem
1.4 (iii) we have $\lambda _n^{(1)}=0$ for all $n\ge 1$. So Theorem
1.5 is true in this case. Now let $f(x_1)>0$. Then $f(x_i)\ge
f(x_1)>0$ for all $i\ge 1$ because $f$ is increasing on the sequence
$\{x_i\}_{i=1}^{\infty}$. Since $f(x_i)\le Cx_i$ for all $i\ge 1$,
we have $0<f(x_i)\le Cx_i$ and so ${1\over {f(x_i)}}\ge {1\over
{Cx_i}}$ for all $i\ge 1$. But by condition (ii),
$\sum_{i=1}^{\infty } {1\over {x_i}}=\infty $. Thus we have
$\sum_{i=1}^{\infty }{1\over {f(x_i)}}=\infty $. The
result in this case then follows immediately from Theorem 1.4 (iii). \ \ \ $\Box$ \\

\noindent{\bf Definition.} ([Hon-Lo]) {\it Let $e$ and $r$ be
positive integers. Let $X=\{x_1, ..., x_e\}$ and $Y=\{y_1, ...,
y_r\}$ be two sets of distinct positive integers. Then we define the
{\it tensor product (set)} of $X$ and $Y$, denoted by $X\odot Y$,
by}
$$
X\odot Y:=\{x_1y_1, ..., x_1y_r, x_2y_1, ..., x_2y_r, ..., x_ey_1,
..., x_ey_r\}.
$$

\noindent{\bf Lemma 2.4.} {\it Let $f$ be a multiplicative function.
Let $e$ and $r$ be positive integers. Let $X=\{x_1, ..., x_e\}$ be a
set of $e$ distinct positive integers such that for any $1\le i\ne
j\le e, \ (x_i, x_j)=1$. Let $Y=\{y_1, ..., y_r\}$ be a set of $r$
distinct positive integers such that for any $1\le i\ne j\le r, \
(y_i, y_j)=1$. Assume that for all $1\le i\le e, 1\le j\le r,\ (x_i,
y_j)=1$. Then the following equality holds:}
$$
(f(X\odot Y))=(f(X))\otimes (f(Y)).
$$

\noindent{\bf Proof.} Since $f$ is multiplicative, we have $f(1)=0$
or $f(1)=1$. If $f(1)=0$, then $f(z)=0$ for every integer $z\ge 1$
because $f$ is multiplicative. Hence we have $ (f(X\odot
Y))=(f(X))\otimes (f(Y))=O_{er}$, the $er\times er$ zero matrix. So
the result holds in this case. Assume now that $f(1)=1$. Then we
have
$$(f(X))=\left(\begin{array}{cccc}
  f(x_1) & 1& ... & 1 \\
  1 & f(x_2)& ... & 1 \\
.. & ..& ..&..\\
 1&1&...& f(x_e)
\end{array}\right)
$$
and
$$
(f(Y))=\left( \begin{array}{cccc}
  f(y_1) & 1& ... & 1 \\
  1 & f(y_2)& ... & 1 \\
.. & ..& ..&..\\
 1&1&...&f(y_r)
\end{array}\right).
$$

Since $f$ is multiplicative, we deduce that
$$
f(x_{i_1}y_{j_1}, \ x_{i_2}y_{j_2})= \left\{\begin{array}{cc}
  f(x_{i_1})f(y_{j_1}) & \ {\rm if} \ i_1=i_2 \ {\rm and}\ j_1=j_2, \\
  f(y_{j_1}) & \ {\rm if} \ i_1\ne i_2 \ {\rm and}\ j_1=j_2, \\
f(x_{i_1}) & \ {\rm if} \ i_1=i_2 \ {\rm and}\ j_1\ne j_2, \\
 1 & \ {\rm if} \ i_1\ne i_2 \ {\rm and}\ j_1\ne j_2.
\end{array}\right.
$$
Thus letting $Y_f=(f(Y))$ gives

$$\begin{array}{rl}
(f(X\odot Y))=& \left(\begin{array}{cccc}
  f(x_1)Y_f & Y_f& ... & Y_f \\
  Y_f & f(x_2)Y_f& ... & Y_f \\
.. & ..& ..&..\\
 Y_f&Y_f&...& f(x_e)Y_f
\end{array}\right)\\
=& \left(\begin{array}{cccc}
  f(x_1) & 1& ... & 1 \\
  1 & f(x_2)& ... & 1 \\
.. & ..& ..&..\\
 1&1&...& f(x_e)
\end{array}\right)\otimes Y_f\\
=& (f(X))\otimes (f(Y))
\end{array}
$$
as required. \ \ \ $\Box$ \\

\noindent{\bf Remark.} If $f$ is {\sc not} multiplicative, then
Lemma 2.2 may fail to be true. For instance, let $X=\{1, 2\}$ and
$Y=\{3, 5\}$. Then $X\odot Y=\{3, 5, 6, 10\}$. Let $f$ be the
arithmetical function defined by $f(l)=l$ for $l\ne 10$ and
$f(10)=9$. Then $f$ is not multiplicative since $f(10)\ne f(2)f(5)$.
On the other hand, we have $((f(X))\otimes (f(Y)))_{44}=10$ and
$(f(X\odot Y))_{44}=9$. This implies that
$(f(X\odot Y))\ne (f(X))\otimes (f(Y)).$\\

We are now in a position to prove Theorem 1.6.\\

\noindent{\bf Proof of Theorem 1.6.} First we prove part (i). For
convenience we let $h:=f_1^{(l_1)}*...*f_c^{(l_c)}*\mu ^{(d)}$ and
$h(p_i(b))=0$ for some $i\ge 1$. Then $(p_i(b), a+be)=1$ or
$(p_i(b), a+b(e+1))=1$. Otherwise we have $p_i(b)|(a+be)$ and
$p_i(b)|(a+b(e+1))$. It implies $p_i(b)|b$ and so $p_i(b)\le b$.
This is absurd since $p_i(b)\ge 1+b$. We may let $(p_i(b),
a+b(e+j))=1$, where $j=0$, or $1$. For any integer $m\ge q$, let
$\nu _m^{(1)}(l_1, ..., l_c, d)\le ...\le \nu _m^{(m)}(l_1, ...,
l_c, d)$ be the eigenvalues of the $m\times m$ matrix $(h(V_m))$
defined on the set
$$V_m:=\{(a+b(e+j))p_i(b), (a+b(e+j))p_i(b)p_{i+w+1}(b)), ...,
(a+b(e+j))p_i(b)p_{i+w+m-1}(b))\},$$ where $w\ge 0$ and
$a+b(e+j)<p_{i+w+1}(b))<...<p_{i+w+m-1}(b)$. Clearly $h$ is
multiplicative since $f_1, ..., f_c$ and $\mu$ are multiplicative.
For each $1\le l\le m-1$, since $p_i(b), a+b(e+j)$ and
$p_{i+w+l}(b)$ are mutually coprime, and note also that
$h(p_i(b))=0$, we have
$h((a+b(e+j))p_i(b)p_{i+w+l}(b))=h(a+b(e+j))h(p_i(b))h(p_{i+w+l}(b))=0$.
Thus we have $(h(V))=O_{m\times m}$, the $m\times m$ zero matrix.
So we have $\nu _m^{(i)}(l_1, ..., l_c, d)=0$ for all $1\le i\le
m$. But by Cauchy's interlacing inequalities we have for any large
enough $n$,
$$\lambda _n^{(q)}(l_1,
..., l_c, d)\le \nu _m^{(q)}(l_1, ..., l_c, d).$$ On the other
hand, Theorem 1.2 (ii) gives $\lambda _n^{(q)}(l_1, ..., l_c,
d)\ge 0.$ So we have $\lambda _n^{(q)}(l_1, ..., l_c, d)=0.$ This
completes the proof of part (i) of Theorem 1.6.

From now on we assume that $h(p_i(b))\ne 0$ for all $i\ge 1$. Next
we prove Theorem 1.6 (ii) for the case $l_1=c=1$ and $d=0$. Then
we have $h=f_1$.

Let $\{1+bt_i\}_{i=0}^{\infty }$ be the sequence consisting of all
those elements in the sequence $\{1+bi\}_{i=0}^{\infty }$ which
are coprime to $a+be$. So $(1+bt_i, a+be)=1$ for all $i\ge 0$.
Then this is an infinite sequence because it contains the set of
all primes strictly greater than $a+be$ in $\{1+bi\}_{i=1}^{\infty
}$, which is infinite by Dirichlet's theorem. For the arithmetic
progression $\{a+bi\}_{i=e}^{\infty }$, consider its subsequence
$$\{a+b(e+(a+be)t_i)\}_{i=0}^{\infty
}=\{(a+be)(1+bt_i)\}_{i=0}^{\infty }.$$ For any integer $m\ge 1$,
let $\gamma _m^{(1)}\le ...\le \gamma _m^{(m)}$ be the eigenvalues
of the $m\times m$ matrix $(f_1(W_m))$ defined on the set
$$W_m:=\{a+be, (a+be)(1+bt_1), ..., (a+be)(1+bt_{m-1})\}$$
and let $\tilde \gamma _m^{(1)}\le ...\le \tilde \gamma _m^{(m)}$
be the eigenvalues of the $m\times m$ matrix $(f_1(\tilde W_m))$
defined on the set
$$\tilde W_m:=\{1, 1+bt_1, ..., 1+bt_{m-1}\}.$$
Since $f_1$ is multiplicative and $(a+be, 1+bt_i)=1$, we have
$(f_1(W_m))=f_1(a+be)(f_1(\tilde W_m))$. So we have $\gamma
_m^{(i)}=f_1(a+be)\tilde \gamma _m^{(i)}$ for $1\le i\le m$. In
particular,
$$ \gamma _m^{(q)}=f_1(a+be)\tilde \gamma _m^{(q)}. \eqno (2-13)$$

Now let $m_n$ be the largest integer $l$ such that
$$t_{l-1}\le \lfloor {{n-1}\over {a+be} }\rfloor ,$$ where $\lfloor
x\rfloor $ denotes the largest integer $\le x$. Clearly
$m_n\rightarrow\infty $ as $n\rightarrow\infty $. Choose $n$ so
that $m_n\ge q$.

By Cauchy's interlacing inequalities
$$
\lambda _n^{(q)}(1, 0)\le \gamma _{m_n}^{(q)}, \eqno (2-14)
$$
and by (2-13) and (2-14),
$$
\lambda _n^{(q)}(1, 0)\le f_1(a+be)\tilde \gamma _{m_n}^{(q)}.
\eqno (2-15)
$$
We claim that ${\rm lim}_{m\rightarrow \infty }\tilde \gamma
_m^{(q)}=0$. Then we have ${\rm lim}_{n\rightarrow \infty }\tilde
\gamma _{m_n}^{(q)}=0$. Thus by Theorem 1.3 and (2-15) we get
${\rm lim}_{n\rightarrow \infty }\lambda _n^{(q)}(1, 0)=0$ as
desired. It remains to prove the assertion which will be done in
the following.

Let $p_1<p_2<...$ denote the primes in the sequence
$\{1+bt_i\}_{i=0}^{\infty }$. Then $\{p_i(b)\}_{i=s}^{\infty
}\subset \{p_i\}_{i=1}^{\infty }$, where $p_{s-1}(b)\le
a+be<p_s(b)$, $s\ge 1$ is an integer and $p_0(b):=1$. Let
$$Q:=\{p_i(b)\}_{i=1}^{\infty}\setminus \{p_i\}_{i=1}^{\infty }.$$
Then $Q$ is a finite set. Since $f_1(p_i(b))\ne 0$ for all $i\ge
1$, we have $\sum_{p\in Q}{1\over {f_1(p)}}<\infty$. So by the
assumption $\sum_{i=1}^{\infty }{1\over {f_1(p_i(b))}}=\infty $ we
have
$$
\sum_{i=1}^{\infty }{1\over {f_1(p_i)}}=\infty . \eqno (2-16)
$$

For $i\ge 1$, let $\pi _i=p_{q-1+i}$. Then $p_{q-1}<\pi _1<...$.
Since $q$ is a fixed number, it follows from (2-16) that
$$
\sum_{i=1}^{\infty }{1\over {f_1(\pi _i)}}=\infty . \eqno (2-17)
$$

Now let $r\ge 2$ be an arbitrary integer and let
$$P_q:=\{1, p_1, ..., p_{q-1}\}, \ \ T_r:=\{1, \pi _1, ..., \pi
_{r-1}\}.$$ It is clear that the matrices $(f_1(P_q))$ and
$(f_1(T_r))$ are positive semi-definite. Consider the tensor
product set $P_q\odot T_r$. Note that the entries in the set
$P_q\odot T_r$ are not arranged in increasing order, but the
eigenvalues of the corresponding matrix do not depend on
rearranging those entries. Since $f_1$ is multiplicative, by Lemma
2.4 we have
$$
(f_1(P_q\odot T_r))=(f_1(P_q))\otimes (f_1(T_r)).
$$

Let $\delta _q^{(1)}\le ...\le \delta _q^{(q)}$ and $\tilde
\lambda _r^{(1)}\le ...\le \tilde \lambda _r^{(r)}$ be the
eigenvalues of the matrix $(f_1(P_q))$ defined on the set $P_q$
and the matrix $(f_1(T_r))$ defined on the set $T_r$ respectively.
Then it is known (see [Hor-J2]) that the eigenvalues of the tensor
product matrix $(f_1(P_q))\otimes (f_1(T_r))$ are given by the set
$$\left\{\delta _q^{(i)}\cdot \tilde \lambda _r^{(j)}\right\}_{1\le j\le
r}^{1\le i\le q}.$$ Notice that
$$
\delta _q^{(1)}\cdot \tilde \lambda _r^{(1)}\le ...\le \delta
_q^{(q)}\cdot \tilde \lambda _r^{(1)}.   \eqno (2-18)
$$

Clearly the sequence $\{1+bt_i\}_{i=0}^{\infty }$ is closed under
the usual multiplication. So the tensor product set $P_q\odot
T_r\subset \{1+bt_i\}_{i=0}^{\infty }$. For any integer $r\ge 2$,
define an integer $m_r$ by
$$m_r:={{p_{q-1}\cdot \pi _{r-1}-1}\over b}+1.$$
Then $P_q\odot T_r\subseteq \{1+bt_i\}_{i=0}^{m_r-1}$. Thus the
matrix $(f_1(P_q\odot T_r))$ defined on $P_q\odot T_r$ is a
principal submatrix of the $m_r\times m_r$ matrix $(f_1(1+bt_i,
1+bt_j))$ defined on the set $\{1, 1+bt_1, ..., 1+bt_{m_r-1}\}$.
Let $\bar \lambda _{qr}^{(1)}\le ...\le \bar \lambda _{qr}^{(qr)}$
be the eigenvalues of $(f_1(P_q\odot T_r))$. Then by Cauchy's
interlacing inequalities we have
$$
\tilde\gamma _{m_r}^{(q)}\le \bar \lambda _{qr}^{(q)}. \eqno
(2-19)
$$
But by (2-18)
$$
\bar \lambda _{qr}^{(q)}\le \delta _q^{(q)}\cdot \tilde \lambda
_r^{(1)}. \eqno (2-20)
$$
So it follows from (2-19) and (2-20) that
$$
\tilde\gamma _{m_r}^{(q)}\le \delta _q^{(q)}\cdot \tilde \lambda
_r^{(1)}. \eqno (2-21)
$$

On the other hand, in Theorem 1.4, if we choose $x=x_1=1$ and
$x_i=\pi _{i-1}$ for $i\ge 2$, then by (2-17) the conditions of
Theorem 1.4 are satisfied. It then follows immediately from
Theorem 1.4 that
$$
{\rm lim}_{r\rightarrow \infty }\tilde \lambda _r^{(1)}=0. \eqno
(2-22)
$$
But by Theorem 1.3 we have that the subsequence $\{\tilde\gamma
_{m_r}^{(q)}\}_{r=1}^{\infty }$ of the sequence $\{\tilde\gamma
_m^{(q)}\}_{m=q}^{\infty }$ converges and
$$
{\rm lim}_{r\rightarrow \infty }\tilde\gamma _{m_r}^{(q)}\ge 0.
\eqno (2-23)
$$
Hence by (2-21)-(2-23), ${\rm lim}_{r\rightarrow \infty
}\tilde\gamma _{m_r}^{(q)}=0$. Finally, again by Theorem 1.2, the
desired result ${\rm lim}_{m\rightarrow \infty }\tilde\gamma
_m^{(q)}=0$ follows immediately. The claim is proved and this
completes the proof of part (ii) of Theorem 1.6 for the case
$l_1=c=1$ and $d=0$.

Now consider part (ii) for the general case. In the case $l_1=c=1$
and $d=0$, we replace $f_1$ by $h=f_1^{(l_1)}*...*f_c^{(l_c)}*\mu
^{(d)}$. Since $f_1, ..., f_c\in {\cal C}$ and $l_1+...+l_c>d$, by
Theorem 1.1 (i) we have $h\in {\cal C}$. Thus $h(p_i(b))\ge 0$ for
all $i\ge 1$. So by the assumption $h(p_i(b))\ne 0$ for all $i\ge
1$ we have that $h(p_i(b))>0$ for all $i\ge 1$. Note that $h$ is
multiplicative. On the other hand, $h$ is increasing on the
sequence $\{{p_i(b)}\}_{i=1}^{\infty}$ because of the formula in
Theorem 1.1 (ii). It remains to prove that
$$
\sum_{i=1}^{\infty }{1\over {h(p_i(b))}}=\infty . \eqno (2-24)
$$
But Theorem 1.1 (ii) tells us
$$
h(p_i(b))=\sum_{j=1}^cl_jf_j(p_i(b))-d<\sum_{j=1}^cl_jf_j(p_i(b)),
$$
So we have
$$
\sum_{i=1}^{\infty }{1\over {h(p_i(b))}}\ge \sum_{i=1}^{\infty
}{1\over {\sum_{j=1}^cl_jf_j(p_i(b))}}\ge {1\over
l}\sum_{i=1}^{\infty }{1\over {\sum_{j=1}^cf_j(p_i(b))}}, \eqno
(2-25)
$$
where $l={\rm max}_{1\le j\le c}l_j$. Hence (2-24) follows
immediately from (2-25) and the condition of Theorem  1.6 (ii). So
Theorem 1.6 (ii) for the general case follows from Theorem 1.6
(ii) for the case $l_1=c=1$ and $d=0$. The proof of Theorem 1.6
(ii) is complete.

Finally Theorem 1.6 (iii) follows from parts (i) and (ii) and
Mertens' Theorem ([Me]). This concludes the proof of Theorem 1.6.
\ \ \ $\Box$ \\

\noindent{\bf 3. Examples}\\

In the present section we give several examples
to demonstrate our main results.\\

\noindent{\bf Example 3.1.} Let $f=\xi _{\varepsilon }$, where $\xi
_{\varepsilon }$ is defined as in the introduction and $\varepsilon
$ is a real number. Then $\xi_{\varepsilon }$ is increasing on any
strictly increasing infinite sequence, and completely multiplicative
if $\varepsilon \ge 0$. Let $J_{\varepsilon }:=\xi_{\varepsilon
}*\mu $. Then $J_{\varepsilon}(1)=1$ and for any integer $m>1$,
$$J_{\varepsilon}(m)=m^{\varepsilon }\prod_{p|m}(1-{1\over
{p^{\varepsilon }}})\ge 0$$ if $\varepsilon\ge 0.$ Thus $\xi
_{\varepsilon }\in {\cal C}_S$ for any set $S$ of positive
integers and so $\xi _{\varepsilon }\in {\cal C}$ for any given
strictly increasing infinite sequence $\{x_i\}_{i=1}^{\infty}$ of
positive integers if $\varepsilon\ge 0$. For integers $c>d\ge 0$,
let ${\lambda}_n^{(1)}(c,d)\le ...\le \lambda _n^{(n)}(c,d)$ be
the eigenvalues of the $n\times n$ matrix $((\xi _{\varepsilon
}^{(c)}*\mu^{(d)})(x_i, x_j))$ defined on the set $S_n=\{x_1, ...,
x_n\}$.

(i). By Theorem 1.4 (ii) we get: If $\varepsilon >0$ and $S_n$
satisfies that for every $1\le i\ne j\le n,\ \ (x_i, x_j)=x$, then
we have
$$x_1^{\varepsilon }-x^{\varepsilon }\le {\lambda}_n^{(1)}(1,0)<
x_1^{\varepsilon }-
x^{\varepsilon }+{{x^{\varepsilon }}\over
{1+\sum_{i=2}^{n}{{x^{\varepsilon }}\over {x_i^{\varepsilon
}-x_1^{\varepsilon }}}}};$$

(ii). ([Hon-Lo]) By Theorem 1.4 (iii) we get: For any given strictly
increasing infinite sequence $\{x_i\}_{i=1}^{\infty }$ consisting of
all but finitely many primes, we have $(x_i, x_j)=1$ for every $i\ne
j$, and by Mertens' theorem ([Me]) we have
$\sum_{i=1}^{\infty}{1\over {x_i^{\varepsilon }}}={\infty}$ if
$\varepsilon \le 1$. So if $0\le \varepsilon \le 1$, then we have
${\rm lim}_{n\rightarrow\infty
}{\lambda}_n^{(1)}(1,0)=x_1^{\varepsilon }-1$;

(iii). By Theorem 1.8 we get: For any given strictly increasing
infinite sequence $\{x_i\}_{i=1}^{\infty }$ of positive integers
which contains the arithmetic progression $\{a+bi\}_{i=e}^{\infty
}$ as its subsequence, where $a, \ b\ge 1$ and $e\ge 0$ are
integers, if $0\le \varepsilon \le 1$, then for any given integer
$q\ge 1$, we have ${\rm
lim}_{n\rightarrow\infty } {\lambda}_n^{(q)}(c,d)=0$.\\

\noindent{\bf Example 3.2.} Let $f=J_{\varepsilon }$, where
$\varepsilon$ is a real number and $J_{\varepsilon}$ is defined in
Example 3.1. Note that if $\varepsilon $ is a positive integer, then
$J_{\varepsilon }$ becomes Jordan's totient function (see, for
example, [A1], [Mc1] or [Mu]). Clearly $J_{\varepsilon }*\mu $ is
multiplicative and $(J_{\varepsilon }*\mu )(1)=1$. It is easy to see
that if $p$ is an odd prime number and $\varepsilon \ge {{{\rm log
}2}\over {{\rm log} 3} }$, then $(J_{\varepsilon }*\mu
)(p)=p^{\varepsilon }-2\ge 0$. For any prime $p$ and integer $l\ge
2$, we have $(J_{\varepsilon }*\mu )(p^l)=p^{(l-2)\varepsilon
}(p^{\varepsilon }-1)^2>0$. Thus $J_{\varepsilon }\in {\cal C}_S$
for any set $S$ of positive odd numbers and so $J_{\varepsilon }\in
{\cal C}$ for any given strictly increasing infinite sequence
$\{x_i\}_{i=1}^{\infty}$ of positive odd numbers if $\varepsilon \ge
{{{\rm log }2}\over {{\rm log} 3} }$. On the other hand, if
$\varepsilon \ge 0$, then for any primes $3\le p_1<p_2$, we have
$J_{\varepsilon }(p_1)=p_1^{\varepsilon }-1\le p_2^{\varepsilon
}-1=J_{\varepsilon }(p_2)$ and for any integer $m\ge 2$, we have $
J_{\varepsilon }(m)\le m^{\varepsilon }$. For integers $c>d\ge 0$,
let ${\lambda}_n^{(1)}(c,d)\le ...\le \lambda _n^{(n)}(c,d)$ be the
eigenvalues of the $n\times n$ matrix $((J_{\varepsilon}^{(c)}*\mu
^{(d)})(x_i, x_j))$ defined on the set $S_n=\{x_1, ..., x_n\}$.

(i). By Theorem 1.4 (ii) we get: For any given strictly increasing
infinite sequence $\{x_i\}_{i=1}^{\infty }$ consisting of all but
finitely many odd primes, if ${{{\rm log}2}\over {{\rm
log}x_1}}\le \varepsilon <1$, then we have
$$x_1^{\varepsilon }-2<{\lambda}_n^{(1)}(1,0)<x_1^{\varepsilon }-2+{1\over
{1+\sum_{i=2}^{n}{1\over {x_i^{\varepsilon }-x_1^{\varepsilon
}}}}};$$ Furthermore by Theorem 1.5, ${\rm
lim}_{n\rightarrow\infty }{\lambda}_n^{(1)}(1,0)=x_1^{\varepsilon
}-2$;

(ii). By Theorem 1.8 we get: For any given strictly increasing
infinite sequence $\{x_i\}_{i=1}^{\infty }$ of positive odd
numbers which contains the arithmetic progression
$\{a+bi\}_{i=e}^{\infty }$ as its subsequence, where $a, \ b\ge 1$
and $e\ge 0$ are integers, if ${{{\rm log}2}\over {{\rm log}3}}\le
\varepsilon <1$, then for any given integer $q\ge 1$, we
have ${\rm lim}_{n\rightarrow\infty }{\lambda}_n^{(q)}(c,d)=0$.\\

\noindent{\bf Example 3.3.} Let $f=\sigma _{\varepsilon }:=\xi
_{\varepsilon }*\xi _0$, where $\varepsilon $ is a real number. Then
for any positive integer $m$ we have
$$\sigma _{\varepsilon }(m)=\sum_{d|m}d^{\varepsilon }.$$
The function $d(m)=\sigma _0(m)$ is the usual divisor function.
The function $\sigma (m)=\sigma _1(m)$ gives the sum of the
divisors of $m$. Clearly $\sigma _{\varepsilon }$ is
multiplicative. Since $\sigma _{\varepsilon }*\mu =\xi
_{\varepsilon }*\xi _0*\mu =\xi _{\varepsilon }$, we have $(\sigma
_{\varepsilon }*\mu )(m)=m^{\varepsilon }>0$ for any integer $m\ge
1$. So $\sigma _{\varepsilon }\in {\cal C}$ for any given strictly
increasing infinite sequence $\{x_i\}_{i=1}^{\infty}$ of positive
integers. Obviously if $\varepsilon \ge 0$ and $p_1<p_2$ are
primes, then $\sigma _{\varepsilon }(p_1)=1+p_1^{\varepsilon }\le
1+p_2^{\varepsilon }=\sigma _{\varepsilon }(p_2)$. For integers
$c>d\ge 0$, let ${\lambda}_n^{(1)}(c,d)\le ...\le \lambda
_n^{(n)}(c,d)$ be the eigenvalues of the $n\times n$ matrix
$((\sigma _{\varepsilon}^{(c)}*\mu ^{(d)})(x_i, x_j))$ defined on
the set $S_n=\{x_1, ..., x_n\}$.

(i). By Theorems 1.4 (ii) we get: For any given strictly
increasing infinite sequence $\{x_i\}_{i=1}^{\infty }$ consisting
of all the primes in ${\bf Z}^+$ except finitely many of them, if
$\varepsilon >0$, then we have
$$x_1^{\varepsilon }<{\lambda}_n^{(1)}(1,0)<x_1^{\varepsilon }+{1\over
{1+\sum_{i=2}^{n}{1\over {x_i^{\varepsilon }-x_1^{\varepsilon
}}}}}.$$ Furthermore by Theorem 1.4 (iii), if $0\le \varepsilon
\le 1$, then we have ${\rm lim}_{n\rightarrow\infty
}{\lambda}_n^{(1)}(1,0)=x_1^{\varepsilon }$;

(ii). By Theorem 1.8 we get: For any given strictly increasing
infinite sequence $\{x_i\}_{i=1}^{\infty }$ of positive integers
which contains the arithmetic progression $\{a+bi\}_{i=e}^{\infty
}$ as its subsequence, where $a, \ b\ge 1$ and $e\ge 0$ are
integers, since
$$\sum_{i=1}^{\infty }{1\over {\sigma _{\varepsilon }
(p_i(b))}}=\sum_{i=1}^{\infty }{1\over {p_i(b)^{\varepsilon
}+1}}\ge {1\over 2}\sum_{i=1}^{\infty }{1\over
{p_i(b)^{\varepsilon }}}$$ if $\varepsilon \ge 0$, we deduce that
if $0\le \varepsilon \le 1$
$$\sum_{i=1}^{\infty }{1\over {\sigma _{\varepsilon }
(p_i(b))}}=\infty .$$ Then for any given integer $q\ge 1$, if
$0\le \varepsilon \le 1$, we have ${\rm lim}_{n\rightarrow\infty
}{\lambda}_n^{(q)}(c,d)=0$.
\\

\noindent{\bf Example 3.4.} Let $f=\psi _{\varepsilon }$, where
$\varepsilon $ is a real number and $\psi _{\varepsilon }$ is
defined for any positive integer $m$ by
$$\psi _{\varepsilon
}(m):=\sum_{d|m}d^{\varepsilon }|\mu ({m\over d})|.$$ The function
$\psi _1$ is called Dedekind's function (see, for instance,
[Mc1]). Clearly $\psi _{\varepsilon }$ is multiplicative. Then for
any positive integer $m$ we have
$$\psi _{\varepsilon
}(m)=m^{\varepsilon }\prod_{p|m}(1+{1\over {p^{\varepsilon
}}})={{J_{2\varepsilon }(m)}\over {J_{\varepsilon }(m)}}.$$
Thus
for any positive integer $l$ and any prime $p$, we have
$$
(\psi _{\varepsilon }*\mu
)(p^l)=\left\{\begin{array}{cc}p^{\varepsilon }, & \ \ {\rm if} \
\ l=1;\\
p^{(l-2)\varepsilon }(p^{2\varepsilon }-1), & \ \ {\rm if } \ \
l\ge 2.
\end{array}\right.
$$
If $\varepsilon \ge 0$, then $\psi _{\varepsilon }\in {\cal C}$
for any given strictly increasing infinite sequences
$\{x_i\}_{i=1}^{\infty}$ of positive integers. For integers
$c>d\ge 0$, let ${\lambda}_n^{(1)}(c,d)\le ...\le \lambda
_n^{(n)}(c,d)$ be the eigenvalues of the $n\times n$ matrix
$((\psi _{\varepsilon}^{(c)}*\mu ^{(d)})(x_i, x_j))$ defined on
the set $S_n=\{x_1, ..., x_n\}$.

(i). By Theorem 1.4 (ii) we get: For any given strictly increasing
infinite sequence $\{x_i\}_{i=1}^{\infty }$ consisting of all the
primes in ${\bf Z}^+$ except finitely many of them, if $\varepsilon
>0$, then we have
$$x_1^{\varepsilon }<{\lambda}_n^{(1)}(1,0)<x_1^{\varepsilon }+{1\over
{1+\sum_{i=2}^{n}{1\over {x_i^{\varepsilon }-x_1^{\varepsilon
}}}}}.$$ Furthermore by Theorem 1.4 (iii), if $0\le \varepsilon
\le 1$, then we have ${\rm lim}_{n\rightarrow\infty
}{\lambda}_n^{(1)}(1,0)=x_1^{\varepsilon }$;

(ii). By Theorem 1.8 we get: For any given strictly increasing
infinite sequence $\{x_i\}_{i=1}^{\infty }$ of positive integers
which contains the arithmetic progression $\{a+bi\}_{i=e}^{\infty
}$ as its subsequence, where $a, \ b\ge 1$ and $e\ge 0$ are
integers, in a same way as in Example 3.3, we can check that for
$\varepsilon \ge 0$, $\psi _{\varepsilon }$ is increasing on the
sequence $\{p_i(b)\}_{i=1}^{\infty}$ and if $0\le \varepsilon \le
1$
$$\sum_{i=1}^{\infty }{1\over {\psi _{\varepsilon
}(p_i(b))}}=\infty .$$ Then for any given integer $q\ge 1$, if
$0\le \varepsilon \le 1$, we have ${\rm lim}_{n\rightarrow\infty
}{\lambda}_n^{(q)}(c,d)=0$.\\

\noindent{\bf Example 3.5.} Let $f=\phi $, Euler's totient function.
Clearly $\phi $ and $\phi *\mu $ are multiplicative, and $\phi
(1)=(\phi *\mu )(1)=1$. For any prime $p$ we have $(\phi *\mu
)(p)=\phi (p)-1=p-2\ge 0$, and for any integer $l\ge 2$ we have
$$(\phi *\mu )(p^l)=\displaystyle \sum_{i=1}^l\phi(p^i)\mu
(p^{l-i})=\phi (p^l)-\phi (p^{l-1})=p^{l-2}(p-1)^2>0.
$$
Thus $\phi \in {\cal C}_S$ for any set $S$ of positive integers
and so $\phi \in {\cal C}$ for any given strictly increasing
infinite sequence $\{x_i\}_{i=1}^{\infty}$ of positive integers.
Note that $\phi (p)=p-1\le p$. So for any primes $p_1<p_2$, $\phi
(p_1)<\phi (p_2)$. For integers $c>d\ge 0$, let
${\lambda}_n^{(1)}(c,d)\le ...\le \lambda _n^{(n)}(c,d)$ be the
eigenvalues of the $n\times n$ matrix $((\phi ^{(c)}*\mu
^{(d)})(x_i, x_j))$ defined on the set $S_n=\{x_1, ..., x_n\}$.

(i). By Theorem 1.4 (ii) we get: For any given strictly increasing
infinite sequence $\{x_i\}_{i=1}^{\infty }$ consisting of all the
primes in ${\bf Z}^+$ except finitely many of them, we have
$$x_1-2<{\lambda}_n^{(1)}(1,0)<x_1-2+{1\over
{1+\sum_{i=2}^{n}{1\over {x_i-x_1}}}}.$$ Furthermore, by Theorem
1.5 we have ${\rm lim}_{n\rightarrow\infty
}{\lambda}_n^{(1)}(1,0)=x_1-2$;

(ii). By Theorem 1.8 we get: For any given strictly increasing
infinite sequence $\{x_i\}_{i=1}^{\infty }$ of positive integers
which contains the arithmetic progression $\{a+bi\}_{i=e}^{\infty
}$ as its subsequence, where $a, \ b\ge 1$ and $e\ge 0$ are
integers, and for any given integer $q\ge 1$, we have
${\rm lim}_{n\rightarrow\infty }{\lambda}_n^{(q)}(c,d)=0$.\\

\noindent{\bf Example 3.6.} Let $f_1=\xi _{\varepsilon }$ and
$f_2=\phi $ be as in Examples 3.1 and 3.5 respectively. Clearly $\xi
_{\varepsilon }$ and $\phi $ are distinct and multiplicative. Note
that $\xi_{\varepsilon }$ is increasing on any strictly increasing
infinite sequence of positive integers if $\varepsilon \ge 0$ and
$\phi $ is increasing on any subsequence of strictly increasing
infinite sequence consisting of all the primes in ${\bf Z}^+$. By
Examples 3.1 and 3.5 we know that $\xi _{\varepsilon }\in {\cal C}$
and $\phi \in {\cal C}$ for any given strictly increasing infinite
sequence $\{x_i\}_{i=1}^{\infty}$ of positive integers if
$\varepsilon\ge 0$. For integers $c_1>0, c_2>0$ and $d\ge 0$, let
${\lambda}_n^{(1)}(c_1,c_2,d)\le ...\le \lambda _n^{(n)}(c_1,c_2,d)$
be the eigenvalues of the $n\times n$ matrix $((\xi _{\varepsilon
}^{(c_1)}*\phi ^{(c_2)}*\mu^{(d)})(x_i, x_j))$ defined on the set
$S_n=\{x_1, ..., x_n\}$. Since for any prime $p$, we have $\phi
(p)\le p$ and $\xi_{\varepsilon}(p)\le p$ if $\varepsilon\le 1$,
then by Theorem 1.7 we get: For any given strictly increasing
infinite sequence $\{x_i\}_{i=1}^{\infty }$ of positive integers
which contains the arithmetic progression $\{a+bi\}_{i=e}^{\infty }$
as its subsequence, where $a, \ b\ge 1$ and $e\ge 0$ are integers,
if $0\le \varepsilon \le 1$ and $c_1+c_2>d$, then for any given
integer $q\ge 1$, we have ${\rm
lim}_{n\rightarrow\infty } {\lambda}_n^{(q)}(c_1,c_2,d)=0$.\\

\noindent{\bf 4. Open questions}\\

Let $\{x_i\}_{i=1}^{\infty}$ be an arbitrary strictly increasing
infinite sequence of positive integers. For an integer $n\ge 1$, let
$S_n=\{x_1,...,x_n\}$. Let $c, q\ge 1$ and $d\ge 0$ be given
integers. Let $\lambda _n^{(1)}(c, d)\le ...\le \lambda _n^{(n)}(c,
d)$ be the eigenvalues of the matrix $((f^{(c)}*\mu ^{(d)})(x_i,
x_j))$ defined on the set $S_n$. It follows from Theorem 1.4 that if
$\{x_i\}_{i=1}^{\infty}$ is a strictly increasing infinite sequence
of positive integers satisfying that for every $i\ne j, \ (x_i,
x_j)=x_1$ and $f\in {\cal C}$ is increasing on the sequence
$\{x_i\}_{i=1}^{\infty}$ and $\sum_{i=1}^{\infty }{1\over
{f(x_i)}}=\infty $, then ${\rm lim}_{n\rightarrow\infty
}{\lambda}_n^{(1)}(1, 0)=0$. Then by Cauchy's interlacing
inequalities and Theorem 1.3 we know that for any given strictly
increasing infinite sequence $\{x_i\}_{i=1}^{\infty }$ of positive
integers which contains a subsequence $\{x'_i\}_{i=1}^{\infty }$
satisfying that for every $i\ne j, \ (x'_i, x'_j)=x'_1$, if $f\in
{\cal C}$ (with respect to the whole sequence $\{x_i\}_{i=1}^{\infty
}$) and $f$ is increasing on the sequence $\{x_i'\}_{i=1}^{\infty}$
and $\sum_{i=1}^{\infty }{1\over {f(x_i')}}=\infty $, then ${\rm
lim}_{n\rightarrow\infty }{\lambda}_n^{(1)}(1, 0)=0$ (Note that this
holds when some $f(x_i')$ is 0). On the other hand, by Theorem 1.8
we know that for any given strictly increasing infinite sequence
$\{x_i\}_{i=1}^{\infty }$ of positive integers containing the
arithmetic progression $\{a+bi\}_{i=e}^{\infty }$ as its
subsequence, if $c>d\ge 0$ and $f\in {\cal C}$ is multiplicative and
increasing on the sequence $\{p_i(b)\}_{i=1}^{\infty}$ and
$\sum_{i=1}^{\infty }{1\over {f(p_i(b))}}=\infty $, where $p_i(b) \
(i\ge 1)$ is defined as in (1-1), then for any given integer $q\ge
1$, we have ${\rm lim}_{n\rightarrow\infty }{\lambda}_n^{(q)}(c,
d)=0$. First we would like to understand for what sequences
$\{x_i\}_{i=1}^{\infty }, \ {\rm lim}_{n\rightarrow\infty
}{\lambda}_n^{(1)}(c, d)=0$.
Namely, we have the following question:\\

\noindent{\bf Question 4.1.} Given any multiplicative function $f$,
and given nonnegative integers $c$ and $d$ such that $c>d$,
characterize all strictly increasing infinite sequences
$\{x_i\}_{i=1}^{\infty }$ of positive integers so that ${\rm
lim}_{n\rightarrow\infty }{\lambda}_n^{(1)}(c, d)=0$, where, as
before, $\lambda _n^{(1)}(c, d)$ is the smallest eigenvalue of the
matrix $((f^{(c)}*\mu ^{(d)})(x_i, x_j))$ defined on the set
$S_n=\{x_1, ..., x_n\}$.\\

Consequently, we raise a further problem.\\

\noindent{\bf Question 4.2.} The same as the previous question, with
${\lambda}_n^{(1)}(c, d)$ is replaced by $\lambda _n^{(q)}(c, d)$,
where, as before, $\lambda _n^{(q)}(c, d)$ is the $q$-th smallest
eigenvalue of the matrix $((f^{(c)}*\mu ^{(d)})(x_i,
x_j))$ defined on the set $S_n=\{x_1, ..., x_n\}$.\\

In concluding this paper we propose the following question and conjecture.\\

\noindent{\bf Question 4.3.} Let $c>d\ge 0$ be given integers and
$\{x_i\}_{i=1}^{\infty}$ be an arbitrary strictly increasing
infinite sequence of positive integers. Let ${\lambda}_n^{(1)}(c,
d)$ be the smallest eigenvalue of the $n\times n$ matrix
$((f^{(c)}*\mu ^{(d)})(x_i, x_j))$ defined on the set $S_n=\{x_1,
..., x_n\}$. Assume that $f\in {\cal C}$ is multiplicative. Are the
following true:

(i). If $f$ satisfies that $f(x_i)\ge Cx_i^{\varepsilon }$ for all
$i\ge 1$, where $\varepsilon >1$ and $C>0$ are constants, do we
have ${\rm lim}_{n\rightarrow\infty }{\lambda}_n^{(1)}(c, d)>0$?

(ii). If $f$ satisfies that $\sum_{i=1}^{\infty }{1\over
{f(x_i)}}<\infty $, do we have
${\rm lim}_{n\rightarrow\infty }{\lambda}_n^{(1)}(c, d)>0$?\\
\\
{\bf Conjecture 4.4.} {\it Let $\varepsilon>0$ and $\{x_i\}^\infty _{i=1}$
be an arbitrary given strictly increasing infinite
sequence of positive integers. Let $\lambda _n^{(1)}$ be the smallest eigenvalue
of the $n\times n$ power GCD matrix $((x_i, x_j)^{\varepsilon})$ defined on the set
$S_n=\{x_1, ..., x_n\}$. If $\sum_{i=1}^\infty \frac{1}{x_i^\varepsilon}<\infty $,
then we have ${\rm lim}_{n\rightarrow\infty }{\lambda}_n^{(1)}>0$.}\\

\noindent{\bf References}\\

\re{[A1]}T.M. Apostol,  {Introduction to analytic number theory},
Springer-Verlag, New York, 1976. \re{[A2]} T.M. Apostol,
Arithmetical properties of generalized Ramanujan sums, {Pacific J.
Math.} {41} (1972), 281-293. \re{[Bh]} R. Bhatia, Infinitely
divisible matrices, {Amer. Math. Monthly} {113} (2006), 221-235.
\re{[Bh-K]} R. Bhatia and H. Kosaki, Mean matrices and infinite
divisibility, {Linear Algebra Appl.} {424} (2007), 36-54.
\re{[Bo-L1]} K. Bourque and S. Ligh, Matrices associated with
classes of arithmetical functions, {J. Number Theory} {45} (1993),
367-376. \re{[Bo-L2]} K. Bourque and S. Ligh, Matrices associated
with arithmetical functions, {Linear Multi-linear Algebra} {34}
(1993), 261-267. \re{[Bo-L3]} K. Bourque and S. Ligh, Matrices
associated with classes of multiplicative functions, {Linear Algebra
Appl.} {216} (1995), 267-275. \re{[Ca]} W. Cao, On Hong's conjecture
for power LCM matrices, {Czechoslovak Math. J.} 57 (2007), 253-268.
\re{[CN]}P. Codec\'a and M. Nair, Calculating a determinant
associated with multiplicative functions, Boll. Unione Mat. Ital.
Sez. B Artic. Ric. Mat. (8) 5 (2002), 545-555. \re{[Co1]} E. Cohen,
A class of arithmetical functions, {Proc. Nat. Acad. Sci. U. S. A.}
{41} (1955), 939-944. \re{[Co2]} E. Cohen, Arithmetical inversion
formulas, Canad. J. Math. 12 (1960), 399-409. \re{[FHZ]} W. Feng, S.
Hong and J. Zhao, Divisibility properties of power LCM matrices by
power GCD matrices on gcd-closed sets, Discrete Math. 309 (2009),
2627-2639. \re{[HK]}P. Haukkanen and I. Korkee, Notes on the
divisibility of LCM and GCD matrices, {International J. Math. and
Math. Science} {6} (2005), 925-935. \re{[He-L-S]} H. Hedenmalm, P.
Lindqvist and K. Seip, A Hilbert space of Dirichlet series and
systems of dilated functions in $L^2(0, 1)$, {Duke Math. J.} {86}
(1997), 1-37. \re{[Hi]} T. Hilberdink, Determinants of
multiplicative Toeplitz matrices, Acta Arith. 125 (2006), 265-284.
\re{[Hon1]} S. Hong, Bounds for determinants of matrices associated
with classes of arithmetical functions, {Linear Algebra Appl.} {281}
(1998), 311-322. \re{[Hon2]} S. Hong, On the Bourque-Ligh conjecture
of least common multiple matrices, {J. Algebra} {218} (1999),
216-228. \re{[Hon3]} S. Hong, Lower bounds for determinants of
matrices associated with classes of arithmetical functions, {Linear
Multilinear Algebra} {45} (1999), 349-358. \re{[Hon4]} S. Hong,
Gcd-closed sets and determinants of matrices associated with
arithmetical functions, {Acta Arith.} {101} (2002),
321-332.\re{[Hon5]} S. Hong, Factorization of matrices associated
with classes of arithmetical functions, {Colloq. Math.} {98} (2003),
113-123. \re{[Hon6]} S. Hong, Notes on power LCM matrices, {Acta
Arith.} {111} (2004), 165-177. \re{[Hon7]} S. Hong, Nonsingularity
of matrices associated with classes of arithmetical functions, {J.
Algebra} {281} (2004), 1-14. \re{[Hon8]} S. Hong, Nonsingularity of
least common multiple matrices on gcd-closed sets, {J. Number
Theory} {113} (2005), 1-9. \re{[Hon9]} S. Hong, Nonsingularity of
matrices associated with classes of arithmetical functions on
lcm-closed sets, {Linear Algebra Appl.} {416} (2006), 124-134.
\re{[Hon10]} S. Hong, Divisibility properties of power GCD matrices
and power LCM matrices, {Linear Algebra Appl.} 428 (2008),
1001-1008. \re{[Hon11]} S. Hong, Asymptotic behavior of the largest
eigenvalue of matrices associated with completely even functions
(mod $r$), Asian-Europ. J. Math. 1 (2008), 225-235. \re{[Hon12]} S.
Hong, Infinite divisibility of Smith matrices, Acta Arith. 134
(2008), 381-386. \re{[Hon-Le]}S. Hong and K.S. Enoch Lee, {
Asymptotic behavior of eigenvalues of reciprocal power LCM
matrices}, Glasgow Math. J. 50 (2008), 163-174. \re{[Hon-Lo]} S.
Hong and R. Loewy, Asymptotic behavior of eigenvalues of greatest
common divisor matrices, {Glasgow Math. J.} {46} (2004), 551-569.
\re{[Hon-S-S]} S. Hong, K.P. Shum and Q. Sun, On nonsingular power
LCM matrices, {Algebra Colloq.} 13 (2006), 689-704. \re{[HonZY]} S.
Hong, J. Zhao and Y. Yin, Divisibility properties of Smith matrices,
Acta Arith. 132 (2008), 161-175. \re{[Hor-J1]} R. Horn and C.R.
Johnson, {Matrix analysis}, Cambridge University Press, 1985.
\re{[Hor-J2]} R. Horn and C.R. Johnson, {Topics in matrix analysis},
Cambridge University Press, 1991. \re{[Hw]} S. Hwang, Cauchy's
interlace theorem for eigenvalues of Hermitian matrices, {Amer.
Math. Monthly} {111} (2004), 157--159. \re{[I-R]} K. Ireland and M.
Rosen, {A classical introduction to modern number theory}, Second
Edition, GTM {84}, Springer-Verlag, New York, 1990. \re{[L]}M. Li,
{Notes on Hong's conjectures of real number power LCM matrices}, J.
Algebra 315 (2007), 654-664. \re{[L-S]} P. Lindqvist and K. Seip,
Note on some greatest common divisor matrices, {Acta Arith.} {84}
(1998), 149-154. \re{[Mc1]} P.J. McCarthy, {Introduction to
arithmetical functions}, Springer-Verlag, New York, 1986. \re{[Mc2]}
P.J. McCarthy, A generalization of Smith's determinant, {Canad.
Math. Bull.} {29} (1986),  109-113. \re{[Me]} F. Mertens, Ein
Beitrag zur analytischen Zahlentheorie, {J. Reine Angew. Math.} {78}
(1874), 46-62. \re{[Mu]} M.R. Murty, {Problems in analytic number
theory}, GTM {206}, Springer-Verlag, New York, 2001. \re{[N]} M.B.
Nathanson, {Elementary methods in number theory,} GTM {195},
Springer-Verlag, 2000. \re{[S]}H.J.S. Smith, On the value of a
certain arithmetical determinant, {Proc. London Math. Soc.} {7}
(1875-1876), 208-212. \re {[T]} Q. Tan, {Divisibility among power
GCD matrices and among power LCM matrices on two coprime divisor
chains}, {Linear Multilinear Algebra} 58 (2010), 659-671.
\re {[TL]}Q. Tan and Z. Lin, Divisibility of determinants of
power GCD matrices and power LCM matrices on finitely many
quasi-coprime divisor chains, Appl. Math. Comput. 217 (2010), 3910-3915.
\re {[TLL]} Q. Tan, Z. Lin and L. Liu, {Divisibility among power GCD
matrices and among power LCM matrices on two coprime divisor chains
II}, {Linear Multilinear Algebra} 59 (2011), 969-983.
\re {[W]} A. Wintner, Diophantine approximations and Hilbert's space,
{Amer. J. Math.} {66} (1944), 564-578.
\re{[XL]} J. Xu and M. Li, Divisibility among power GCD matrices and
among power LCM matrices on three coprime divisor chains, Linear
Multilinear Algebra 59 (2011), 773-788.

\end{document}